\input amstex
\documentstyle {amsppt}
\pageheight{50.5pc}
\pagewidth{32pc}
\def\mbR{\Bbb R}
\def\cF{\Cal F}

\def\wt{\widetilde}

\topmatter
\title{}
{
One version of the Clark representation theorem for Arratia flow
}
\endtitle
\author
Andrey A. Dorogovtsev
\endauthor
\address
Institute of Mathematics,National Academy of
Science of Ukraine, Tereshchenkivska 3, 10601, Kiev, Ukraine
\endaddress
\keywords
Brownian motion, coalescence, Clark representation, stochastic integral
\endkeywords
\subjclass 60H05, 60H40, 60J65, 60K35
\endsubjclass
\email
adoro\@imath.kiev.ua
\endemail
\abstract The article contains description of the functionals from
the family of coalescing Brownian particles. New type of the
stochastic integral is introduced and used.
\endabstract

\endtopmatter
\rightheadtext{ One version of the Clark representation theorem for
Arratia flow } \leftheadtext{ Andrey A. Dorogovtsev }

\document

\head
 Introduction
\endhead

The aim of this article is to establish the Clark representation
for the functionals from the Arratia flow of coalescing Brownian
particles \cite{1-4}. The following description of this flow will
be used. We consider the random process $\{x(u);u\in \Bbb R\}$
with the values in $C([0;1])$ such, that for every $u_1<...<u_n$

1) $x(u_k,\cdot )$ is the standard Wiener process starting at the
point $u_k,$

2) $\forall t\in [0;1]$
$$x(u_1,t )\leq ...\leq x(u_n,t ),$$

3) The distribution of $(x(u_1,\cdot ), ..., x(u_n,\cdot ))$
coincides with the distribution of the standard $n$-dimensional
Wiener process starting at $(u_1,...,u_n)$ on the set
$$ \{f\in C([0;1], \Bbb R^n): f_k(0)=u_k, k=1,...,n, f_1(t)<...<f_n(t),
t\in [0;1]\}.$$

Roughly speaking the process $x$ can be described as a family of
Wiener particles which start from every point of $\Bbb R,$ move
independently up to the moment of the meeting then coalesce and move
together.

The following fact is well-known. If $\{w(t); t\in[0;1]\}$  is a
standard  Wiener process and square-integrable random variable $\alpha$   is
measurable with respect to $w,$  then $\alpha$  can be represented as
a sum
$$
\alpha=E\alpha+\int^1_0f(t)dw(t),
$$
with the usual Ito stochastic integral in the right side. Our aim is
to establish the variant of this theorem for the case, when $\alpha$
is measurable with respect to Arratia flow $\{x(u,t); u\in[0;U],
t\in[0;1]\}.$ It follows from the description above, that the
Brownian motions $\{x(u,\cdot); u\in[0;U]\}$  are not jointly
Gaussian. So, the original Clark theorem can not be used in this
situation. The article is divided onto three parts. In the first
part the construction of the stochastic integral with respect to
Arratia flow is presented. The next part is devoted to the variants
of the Clark theorem for finite number of the Brownian motions
stopped in the random times. In the last part the modification of
the construction from the first part is applied to the
representation of the functionals from the flow. \head 1. Spatial
stochastic integral with respect to Arratia flow
\endhead

Let $\{x(u);u\in \Bbb R\}$ be the Arratia flow, i.e. the flow of
Brownian particles with coalescence described above. For $U
>0$ consider a partition $\pi$ of the interval $[0;U],\ \pi =\{
u_0=0,...,u_n=U\}.$ For $k=1,...,n$ define
$$\tau (u_k)=\inf \{1,t\in [0;1]: x(u_k,t)=x(u_{k-1},t)\}.$$
Note, that $\tau (u_k),\ k=1,...,n$ are stopping moments with
respect to the flow\
$$F_t^\pi = \sigma (x(u_k,s),\ k=1,...,n,\ s\leq t).$$
Let us consider for a bounded measurable function $a :\Bbb
R\rightarrow \Bbb R$ the sum
$$S_\pi=\sum _{k=1}^n \int _0^{\tau (u_k)}a
(x(u_k,s))dx(u_k,s).\tag1.1$$ Our aim is to investigate the limit of
$S_\pi$ under
$$|\pi |=\max _{k=0,...,n-1}(u_{k+1}-u_k)\rightarrow 0$$
and its properties depending on the function $a $ and the spatial
variable $U.$

Let us begin with the moments of $S_\pi.$ It follows from the
standard properties of Ito stochastic integral, that
$$ES_\pi=0,$$
and
$$ES_\pi^2=E\sum _{k=1}^n\int _0^{\tau (u_k)}a ^2
(x(u_k,s))ds.\tag1.2$$ Let us denote

$$ \overline{S}_\pi=\sum _{k=1}^n\int _0^{\tau (u_k)}a ^2
(x(u_k,s))ds.$$

Consider the sequence of increasing partitions $\{\pi _n;n\geq 1\}$
of the interval $[0;U]$ with  $|\pi _n|\rightarrow 0, n\rightarrow
\infty .$

\proclaim {Lemma 1.1} There exists a limit
$$\lim _{n\rightarrow \infty }\overline{S}_{\pi _n}\    a.s.
\tag1.3$$
\endproclaim
\demo{Proof}
 To prove the lemma we will check two properties of the
sequence $\{\overline{S}_{\pi _n};\ n\geq 1\}:$
$$\forall n\geq 1: \overline{S}_{\pi _n}\leq\overline{S}_{\pi
_{n+1}},\tag1.4$$ and
$$\sup _{n\geq 1}E\overline{S}_{\pi _n}<+\infty.\tag1.5$$
Note that it is enough to prove (1.4) in the case, when $\pi _{n+1}$
contains only one additional point $v_0$ comparing with $\pi _n.$
Suppose that $\pi _n =\{ u_0=0,...,u_n=U\}$ and $\pi _{n+1} =\{
u_0=0,...,u_k,v_0,u_{k+1},...,u_n=U\}.$ Denote
$$\acute{\tau} (u_{k+1})=\inf \{1,t\in [0;1]: x(u_{k+1},t)=x(v_0,t)\}. $$
Now
$$\aligned
\overline{S}_{\pi_{n+1}}-\overline{S}_{\pi _n}=&\int
_0^{\acute{\tau} (u_{k+1})}a ^2 (x(u_{k+1},s))ds+\\
&+\int _0^{\tau (v_0)}a ^2 (x(v_0,s))ds-\int _0^{\tau (u_{k+1})}a ^2
(x(u_{k+1},s))ds. \endaligned$$
There are two possibilities. In the
first one $\tau (v_0)<\tau (u_{k+1}).$ Now $\acute{\tau}
(u_{k+1})=\tau (u_{k+1}).$ So in this case
$$\overline{S}_{\pi_{n+1}}-\overline{S}_{\pi _n}=\int _0^{\tau (v_0)}a ^2
(x(v_0,s))ds\geq 0.$$ The next case is $\tau (v_0)\geq \tau
(u_{k+1}).$ This possibility can be realized only if $\tau (v_0)=
\tau (u_{k+1}).$ Now $\acute{\tau} (u_{k+1})\leq \tau (v_0)$ and
$$\int _0^{\tau (u_{k+1})}a ^2 (x(u_{k+1},s))ds=\int
_0^{\acute{\tau} (u_{k+1})}a ^2 (x(u_{k+1},s))ds+\int_{\acute{\tau}
(u_{k+1})}^{\tau (u_{k+1})}a ^2 (x(v_0,s))ds.$$ So, in this case
$$\aligned
\overline{S}_{\pi_{n+1}}&-\overline{S}_{\pi _n}=\\
&=\int _0^{\tau (v_0)}a ^2(x(v_0,s))ds-\int_{\acute{\tau}
(u_{k+1})}^{\tau (u_{k+1})}a ^2 (x(v_0,s))ds=\\
&=\int _0^{\tau (v_0)}a^2(x(v_0,s))ds-\int_{\acute{\tau}
(u_{k+1})}^{\tau (v_0)}a ^2 (x(v_0,s))ds=\\
&=\int _0^{\acute{\tau} (u_{k+1})}a ^2 (x(v_0,s))ds\geq 0.
\endaligned$$ Hence (1.4) is true. Let us
estimate the expectation of $\overline{S}_{\pi _n}.$ Consider two
independent standard Wiener processes $w_1,w_2$ which start from 0
and $u>0$ correspondingly. Denote
$$\tau=\inf \{1,\ t:w_1(t)=w_2(t)\}.$$
Then
$$E\tau=\int _0^1\int _{-u}^up_{2t}(v)dv+\int _{-u}^up_2(v)dv,\tag1.6$$
where $p_t$ is the density of the normal distribution with zero mean
and covariance $t.$ It follows from (1.6) that
$$E\tau\sim \frac{3u}{2\sqrt{\pi}},\ u\rightarrow 0+.\tag1.7$$
Consequently,
$$\overline{\lim _{n\rightarrow \infty}}E\overline{S}_{\pi _n}
\leq U\frac{3}{2\sqrt{\pi}}\sup_\Bbb Ra ^2.$$
Now the statement of the lemma follows from (1.4) and
(1.7).
\enddemo

\remark{Remark 1} It follows from the proof of the lemma that there
exists a limit
$$\lim _{n\rightarrow \infty}E\overline{S}_{\pi_n}.$$
\endremark

\proclaim {Lemma 1.2} There exists a limit
$$m(U)=L_2-\lim _{n\rightarrow \infty}S_{\pi _n}.$$
\endproclaim
\demo {Proof} Let the partitions $\pi _n,\pi _{n+1}$ be the same as
in the proof of the previous lemma. Then
$$\aligned
ES_{\pi _n}S_{\pi _{n+1}}=&E\sum _{j_1=1}^n\int _0^{\tau
(u_{j_1})}a (x(u_{j_1},s))dx(u_{j_1},s)\cdot\\
&\cdot(\sum _{j_2\neq k+1}\int _0^{\tau (u_{j_2})}a
(x(u_{j_2},s))dx(u_{j_2},s)+\int _0^{\tau (v_0)}a
(x(v_0,s))dx(v_0,s)+\\
&+\int _0^{\acute{\tau} (u_{k+1})}a ^2
(x(u_{k+1},s))dx(u_{k+1},s))=\\
=&E\sum _{j\neq k+1}\int _0^{\tau (u_{j})}a ^2(x(u_{j},s))ds+E\int
_0^{\tau (u_{k+1})}a (x(u_{k+1},s))dx(u_{k+1},s)\cdot\\
&\cdot(\int _0^{\tau (v_0)}a (x(v_0,s))dx(v_0,s)+\int
_0^{\acute{\tau} (u_{k+1})}a (x(u_{k+1},s))dx(u_{k+1},s)))=\\
=&E\sum _{j\neq k+1}\int _0^{\tau (u_{j})}a ^2(x(u_{j},s))ds+ E\int
_0^{\tau (u_{k+1})\wedge \acute{\tau} (u_{k+1})}a
^2(x(u_{k+1},s))ds+\\
 &+E\int _{\acute{\tau} (u_{k+1})\wedge \tau (v_0)}^{\tau
(u_{k+1})\wedge \tau (v_0)}a ^2 (x(u_{k+1},s))ds=\\
=&E\sum _{j\neq k+1}\int _0^{\tau (u_{j})}a ^2(x(u_{j},s))ds+ E\int
_0^{\tau (u_{k+1})}a ^2 (x(u_{k+1},s))ds=\\
=&ES^2_{\pi _n}.
\endaligned$$
Consequently for all $n\leq m$
$$ES_{\pi _n}S_{\pi _m}=E\overline{S}_{\pi _n}.$$
Now the statement of the lemma follows from the remark 1.
\enddemo

\remark {Remark 2} Note, that the limit $m(U)$ does not depend on
the choice of the sequence of partitions $\{\pi _n;n\geq 1\}.$
\endremark

To prove this we need in an estimation of the rate of convergence
$ES^2_{\pi _n}$ to its limit.

\proclaim {Lemma 1.3} There exists a constant $C$ such, that for
every partition $\pi$ of the interval $[0;U]$
$$|ES^2_\pi-Em^2(U)|\leq C|\pi|\sup_\Bbb Ra ^2.\tag1.8$$
\endproclaim
\demo{Proof} First consider the partitions $\pi ',\pi ''$ where $\pi
''$ is obtained from $\pi '$ by adding one point on the interval
$[u_k,u_{k+1}].$ As it was mentioned in the proof of the lemma 1
$$\overline{S}_{\pi ''}-\overline{S}_{\pi '}=\int _0^{\zeta (v_0)}a ^2
(x(v_0,s))ds.\tag1.9$$ Here
$$\zeta (v_0)=\inf \{1;t:(x(v_0,t)-x(u_k,t))(x(v_0,t)-x(u_{k+1},t))=0\}.$$
Let us estimate $E\zeta (v_0).$ Consider the standard Wiener process
$\overrightarrow{w}$ on the plane, which is starting from the point
$\overrightarrow{r}$. Suppose that this point lies inside the angle
with the vertex in the origin. Let the value of the angle be less
then $\frac{\pi}{2}$ and the angle lies in the part of the plane
where the both coordinate are nonnegative. Define $\acute{\zeta}$
the first exit time of $\overrightarrow{w}$ from the angle. Then
enlarging the angle up to $\frac{\pi}{2}$ and using one-dimension
expressions like (1.6) we can check that there exists $C> 0$ such,
that
$$E\acute{\zeta}\wedge 1\leq Cr_1r_2,\tag1.10$$
where $\overrightarrow{r}=(r_1,r_2).$

>From this remarks we can conclude that there exists $C_1> 0$ such,
that
$$E\zeta (v_0)\leq C_1(u_{k+1}-v_0)(v_0-u_k).\tag1.11$$
This conclusion can be obtained if we note, that $\zeta (v_0)$ is
the minimum of 1 and the first exit time of the 3-dimensional Wiener
process from the space angle with the value $\frac{\pi}{3}.$ It
follows from (1.9) and (1.11) that
$$E(\overline{S}_{\pi ''}-\overline{S}_{\pi '})\leq C_1
\sup_\Bbb Ra ^2(u_{k+1}-v_0)(v_0-u_k).\tag1.12$$ Now let us consider
the general case when $\pi ''$ is obtained from $\pi '$ by the
adding of a few new points. Denote by $v_1<...<v_m$ the new points
on the interval $[u_k,u_{k+1}].$ Then the new amount which is
obtained in $E(\overline{S}_{\pi ''}-\overline{S}_{\pi '})$ from
this points can be estimated due to (1.12) by the sum
$$C_1\sup_\Bbb Ra ^2\sum _{j=1}^m(v_j-v_{j-1})(u_{k+1}-v_j),$$
where we suppose, that $v_0=u_k.$ Consequently
$$E(\overline{S}_{\pi ''}-\overline{S}_{\pi '})\leq C_1\sup_\Bbb Ra ^2\sum _{k=0}^{n-1}
(u_{k+1}-u_k)^2\leq C_1\sup_\Bbb Ra ^2U|\pi '|.$$ This inequality
leads to the existence of the limit
$$\lim _{|\pi|\rightarrow 0}E\overline{S}_{\pi }.$$
It follows from the proof of the lemma 2, that
$$E m^2(U)=\lim _{|\pi|\rightarrow 0}E\overline{S}_{\pi }.$$
It is clear now that (1.8) holds.
\enddemo

The independence $m(U)$ from a choice of the sequence $\{\pi
_n;n\geq 1\}$ now follows in standard way.

Define for $U\geq 0$ the $\sigma -$field
$$\widetilde{F}_U=\sigma (x(u,\cdot );0\leq u\leq U).$$

\proclaim {Lemma 1.4} The process $\{m(U);U\geq 0\}$ is
$(\widetilde{F}_U)-$martingale.
\endproclaim
\demo{Proof} The measurability of $m(U)$ with respect to
$\widetilde{F}_U$ is evident. Let $0\leq U_1<U_2.$ Consider the
partition $|\pi|$ of $[0;U_2]$ which contains the point $U_1.$ Then
$$\aligned
E(S_\pi /\widetilde{F}_{U_1})&=\sum _{u_k\leq U_1} \int _0^{\tau
(u_k)}a (x(u_k,s))dx(u_k,s)+\\
&+E(\sum _{u_k>U_1} \int _0^{\tau (u_k)}a
(x(u_k,s))dx(u_k,s)/\widetilde{F}_{U_1}).
\endaligned$$ To prove
that the last summand is equal to zero it is enough to consider the
expression
$$E(\int _0^{\tau (u)}a (x(u,s))dx(u,s)/\widetilde{F}_{U_1})$$
for $u>U_1.$ Take $0\leq u_1<...<u_n=U_1.$ For a bounded Borel
function $f:C([0;1])^n\rightarrow \Bbb R$ the expectation
$$E\int _0^{\tau (u)}a (x(u,s))dx(u,s)f(x(u_1,\cdot ),...,x(u_n,\cdot ))\tag1.13$$
can be rewritten as
$$E\int _0^{\acute{\tau}}a (w(s))dw(s)\acute{f}(w_1,...,w_n),$$
where $w$ and $w_1,...,w_n$ are independent standard Wiener
processes starting from the points $u$ and $u_1,...,u_n$
correspondingly, and $\acute{\tau}$ is a stopping time for
$(w,w_1,...,w_n), \acute{f}$ is a bounded Borel function on
$C([0;1])^n.$ Denote by $\Gamma$ the $\sigma -$field corresponding
to $\acute{\tau}.$ Then
$$\aligned
E\int _0^{\acute{\tau}}&a (w(s))dw(s)\acute{f}(w_1,...,w_n)=\\
&=E\int _0^{\acute{\tau}}a
(w(s))dw(s)E(\acute{f}(w_1,...,w_n)/\Gamma )=\\
&=E\int _0^{\acute{\tau}}a
(w(s))dw(s)\tilde{f}(\tilde{w}_1,...,\tilde{w}_n),
\endaligned$$
where $\tilde{w}_k(s)=w_k(s\wedge \acute{\tau}),k=1,...,n,$ and
$\tilde{f}$ is new bounded Borel function. Due to the Clark
representation theorem
$$\tilde{f}(\tilde{w}_1,...,\tilde{w}_n)=c+\sum _{k=1}^n\int _0^{\acute{\tau}}
\eta _k(s)dw_k(s),$$ where for $k=1,...,n   \eta _k$ is the
square-integrable random function adapted to the flow
$$\Gamma _t=\sigma (w(s),w_1(s),...,w_n(s),s\leq t).$$
Consequently,
$$\aligned
E\int _0^{\acute{\tau}}a
(w(s))&dw(s)\tilde{f}(\tilde{w}_1,...,\tilde{w}_n)=E\int
_0^{\acute{\tau}}a (w(s))dw(s)\cdot\\
&\cdot\left(c+\sum _{k=1}^n\int _0^{\acute{\tau}} \eta
_k(s)dw_k(s)\right)=0.
\endaligned$$ Hence (1.13) also equal to
zero. Finally
$$E(\sum _{k>U_1} \int _0^{\tau (u_k)}a
(x(u_k,s))dx(u_k,s)/\widetilde{F}_{U_1})=0.$$ Taking the limit under
the diameter of partition tends to 0 we get the statement of the
lemma.
\enddemo

\head
2. Clark representation for the finite family of coalescing Brownian
motions
\endhead
This section is devoted to the integral representation of the
functionals from $x(u_1,\cdot),\ldots, x(u_n,\cdot),$
$u_1<u_2<\ldots<u_n.$ Let us start with the following simple
lemma, which was already used in the previous section.  \proclaim
{Lemma 2.1} Let $w$ be the standard Wiener process on $[0;1]$  and
$0\leq\tau\leq1$ be the stopping time for $w.$ Suppose, that the
square-integrable random variable $\alpha$ is measurable with
respect to $\{w(\tau\wedge t); t\in[0;1]\}.$ Then $\alpha$  can be
represented as
$$
\alpha=E\alpha+\int^\tau_0f(t)dw(t)
$$
with the certain adapted square-integrable random function $f.$
\endproclaim
\demo{Proof} Note, that
$$
\sigma(w(\tau\wedge t); t\in[0;1])=\cF_\tau,
$$
where $\cF_\tau$   is the $\sigma$-field corresponding to the stopping moment
$\tau.$  Now, due to the original Clark theorem
$$
\alpha=E\alpha+\int^1_0f(t)dw(t).
$$
It remains now to apply the conditional expectation with respect
to $\cF_\tau$ to the both sides of this equality. Lemma is proved.
\enddemo

Consider the following situation. Let $w_1, w_2$  be an independent
standard Wiener processes on $[0;1]$  and $\tau$  be a stopping time with
respect to its join flow of $\sigma$-fields. The processes $w_1, w_2$
and the random variable $\tau$  can be considered on the product of
probability spaces $\Omega_1\times\Omega_2.$  Here $\Omega_1$   is
related to $w_1$  and $\Omega_2$  is related to $w_2.$
\proclaim
{Lemma 2.2}
For every fixed $\omega_1\in\Omega_1$  the random variable
$\tau(\omega_1,\cdot)$  on $\Omega_2$  is the stopping moment for
$w_2$  on $\Omega_2.$
\endproclaim
\demo
{Proof}
The set
$\{\omega_2: \tau(\omega_1,\omega_2)<t\}$   is the cross section of
$\{\tau<t\}$  in $\Omega_1\times\Omega_2.$   Hence its
measurability with respect to $\sigma(w_2(s); s\leq t)$  follows from
the usual arguments of measure theory.
\enddemo

The previous two lemmas lead to the following result.
\proclaim
{Theorem 2.1}
Let $w_0, w_1,\ldots, w_n$   be an independent standard Wiener processes
on $[0;1]$  and for every $k=1,\ldots,n $  \ $\tau_k$  is the stopping time
for the process $(w_0, w_1,\ldots, w_k).$  Suppose, that the
square-integrable random variable $\alpha$  is measurable with
respect to the set
$(w_0(\cdot), w_1(\tau_1\wedge\cdot),\ldots,w_n(\tau_n\wedge\cdot)).$
Then $\alpha$   can be represented as
$$
\alpha=E\alpha+\sum^n_{k=0}\int^{\tau_k}_0f_k(t)dw_k(t),
$$
where $\tau_0=1$   and $f_k$  is adapted to the flow generated by $w_k$
under fixed $w_j, j\ne k.$
\endproclaim
\demo{Proof} Denote $\wt{w}_k(t)=w_k(\tau_k\wedge t), \
k=1,\ldots,n.$  Consider the random variable
$\alpha-E(\alpha/\wt{w}_0,\ldots,\wt{w}_{n-1}).$   It is measurable
with respect to $\wt{w}_n$ under fixed
$\wt{w}_0,\ldots,\wt{w}_{n-1}$ and has zero mean. Due to the
previous lemma it can be written as
$$
\alpha-E(\alpha/\wt{w}_0,\ldots,\wt{w}_{n-1})=\int^{\tau_n}_0f_n(t)
dw_n(t),
$$
where the random function $f_n$  under fixed
$\wt{w}_0,\ldots,\wt{w}_{n-1}$
is adapted to the flow generated by $w_n.$  Repeat the same
procedure to the random variable
$E(\alpha/\wt{w}_0,\ldots,\wt{w}_{n-1}).$   Then
$$
E(\alpha/\wt{w}_0,\ldots,\wt{w}_{n-1})-
E(\alpha/\wt{w}_0,\ldots,\wt{w}_{n-2})=\int^{\tau_{n-1}}_0f_{n-1}(t)dw_{n-1}
(t).
$$
After $n$  steps we will get the statement of the theorem.
\enddemo

\remark{Remark}
 Note, that the representation from the theorem  has the following
property
$$
E\alpha^2=(E\alpha)^2+\sum^n_{k=1}E\int^{\tau_k}_0f_k(t)^2dt.
$$
\endremark

Consider an example of application of the theorem 2.1.
\example
{Example 2.1}   Let $\alpha$ be the square-integrable random variable
measurable with respect to
$x(u_0, \cdot), \ldots, x(u_n, \cdot),$
where $u_0, \ldots, u_n $   are the different points. Define the
random moments
$$
\tau_0=1, \
\tau_k=\inf\{1, t: \ x(u_k,t)\in\{x(u_0,t), \ldots, x(u_{k-1},t)\}\},
 k=1,\ldots, n.
$$
Then $\alpha$   can be represented as
$$
\alpha=E\alpha+\sum^n_{k=0}\int^{\tau_k}_0f_k(t)dx(u_k,t),
$$
where for every $k$  the random function $f_k$  is measurable with respect
$x(u_0,\cdot),\ldots, x(u_k,\cdot)$
and (under fixed
$x(u_0,\cdot),\ldots, x(u_{k-1},\cdot)$)
is adapted to the flow $x(u_k, \tau_k\wedge\cdot).$   In this representation
$$
E\alpha^2=(E\alpha)^2+\sum^n_{k=0}E
\int^{\tau_k}_0f_k(t)^2dt.
$$
\endexample
\head
3. Clark representation
\endhead

Let $\alpha$   be the square-integrable random variable measurable
with respect to $\{x(u,\cdot); u\in[0;U]\}.$   Suppose, that
$\{u_n; n\geq0\}$  is a dense set in $[0;U]$ containing 0 and $U.$
Define the random moments $\{\tau_k; k\geq0\}$   as in the theorem
2.1. The following analog of the Clark representation holds.
\proclaim {Theorem 3.1} The random variable $\alpha$   can be
represented as an infinite sum
$$
\alpha=E\alpha+\sum^\infty_{n=0}\int^{\tau_k}_0f_k(t)dx(u_k,t),
\tag3.1
$$
where $\{f_k\}$   satisfy the same conditions as in the theorem 2.1   and
the series converges in the square mean. Moreover
$$
E\alpha^2=(E\alpha)^2+\sum^\infty_{n=0}E\int^{\tau_k}_0f_k(t)dt.
$$
\endproclaim
\demo {Proof} As it was mentioned in the first section $x$   has a
cadl\'ag trajectories as a random process in $C([0;1]).$
Consequently
$$
\sigma(x(u,\cdot); u\in[0;U])=\sigma(x(u_n;\cdot); n\geq 0)=
\bigvee^\infty_{n=0}\sigma(x(u_0,\cdot),\ldots,x(u_n,\cdot)).
$$
Hence due to the Levy theorem
$$
\alpha=L_2\hbox{-}\lim_{n\to\infty}
E(\alpha/x(u_0,\cdot),\ldots,x(u_n,\cdot)).
$$
Due to the theorem 2.1
$$
E(\alpha/x(u_0,\cdot),\ldots,x(u_n,\cdot))=
E\alpha+\sum^n_{k=0}\int^{\tau_k}_0f_k(t)dx(u_k,t),
$$
where $\tau_k$   and $f_k$  do not change with $n.$   So, taking the
limit under $n\to\infty$  we get the statement of the theorem.
\enddemo

Note, that the sum in (3.1)  is closely related to the spatial
stochastic integral which was built in the first section. Really,
suppose, that $a$ is bounded measurable function on $\mbR$  and
the set $\{u_n; n\geq0\}$  is dense in $[0;U]$   with $u_0=0, \
{u}_1=U.$ \proclaim {Lemma 3.1}
$$
\sum^\infty_{n=0}\int^{\tau_n}_0a(x(u_n,t))dx(u_n,t)=m(U)+\int
_0^1a(x(0,t))dx(0,t),
$$
where $m(U)$  was defined in the first section.
\endproclaim
\demo{Proof} Note, that for every $n\geq1$  the points $u_0,\ldots,
u_n$   if ordered in the growing order form a partition of $[0;U].$
Under $n\to\infty$   these partitions increase and their diameters
tend to zero. To prove the lemma it remains to note that for every
$n\geq1$   the sum $ \sum^n_{k=0}\int^{\tau_k}_0a(x(u_k,t))dx(u_k,t)
$ consider with the sum $S_\pi$  for the corresponding partition.
Lemma is proved.

\enddemo

\enddocument